\documentclass[english]{smfart}
\usepackage{smfenum,bull,graphicx,upref,xspace,amssymb,amscd,array,url,smfthm}
\usepackage[all,dvips,cmtip]{xy}

\title[Gauss-Manin systems and~Frobenius~structures~(II)]
{Gauss-Manin systems, Brieskorn lattices and Frobenius structures (II)}

\alttitle{Syst\`emes de Gauss-Manin, r\'eseaux de Brieskorn et structures de Frobenius (II)}

\author[A.~Douai]{Antoine Douai}
\address{UMR 6621 du CNRS\\ Laboratoire J.A. Dieudonn\'e\\
Universit\'e de Nice\\ Parc Valrose\\ 06108 Nice cedex 2\\ France}
\email{douai@math.unice.fr}

\author[C.~Sabbah]{Claude Sabbah}
\address{UMR 7640 du CNRS\\
Centre de Math{\'e}matiques\\
{\'E}cole polytechnique\\
F--91128 Palaiseau cedex\\
France}
\email{sabbah@math.polytechnique.fr}
\urladdr{http://www.math.polytechnique.fr/cmat/sabbah/sabbah.html}

\makeatletter
\def\mainmatter{\renewcommand{\baselinestretch}{1.1}\normalfont}
\def\backmatter{\renewcommand{\baselinestretch}{1}\normalfont}

\numberwithin{equation}{section}
\NumberTheoremsAs{equation}
\def\thesubsection{\thesection.\alph{subsection}}

\newdir^{ (}{{}*!/-5pt/\dir^{(}}
\newdir_{ (}{{}*!/-5pt/\dir_{(}}

\def\@tempa{english}
\ifx\@tempa\smf@language
   \def\remsname{Remarks}%
   \def\exemsname{Examples}%
   \def\reelname{Re}%
\else
   \def\remsname{Remarques}%
   \def\exemsname{Exemples}%
   \def\reelname{R\acute{e}}%
\fi

\newcommand{\RedefinitSymbole}[1]{%
\expandafter\let\csname old\string#1\endcsname=#1
\let#1=\relax
\newcommand{#1}{\csname old\string#1\endcsname\,}%
}
\RedefinitSymbole{\forall} \RedefinitSymbole{\exists}

\newenvironment{proposition}{\begin{prop}}{\end{prop}}
\newenvironment{lemme}{\begin{lemm}}{\end{lemm}}
\newenvironment{corollaire}{\begin{coro}}{\end{coro}}

\newenvironment{remarque}{\begin{rema}}{\end{rema}}

\newenvironment{exemple}{\begin{exem}}{\end{exem}}

\newenvironment{Proof}{\skippointrait\begin{proof}}{\end{proof}}

\theoremstyle{plain}
\newtheorem*{assertion*}{Assertion}
\newtheorem{mtheoreme}{\theoname}

\theoremstyle{remark}
\newtheorem*{remarque*}{\remaname}

\let\wh\widehat
\let\wt\widetilde
\let\ov\overline
\let\emptyset\varnothing
\let\ldots\dots

\let\moins\smallsetminus
\let\epsilon\varepsilon

\let\leq\leqslant
\let\geq\geqslant

\newcommand{\bbullet}{{\scriptscriptstyle\bullet}}

\newcommand{\ooplus}{\mathop\oplus\limits}
\newcommand{\ootimes}{\mathop\otimes\limits}

\newcommand{\module}[1]{\left\vert#1\right\vert}

\DeclareMathOperator{\id}{Id}

\newcommand{\isom}{\stackrel{\sim}{\longrightarrow}}

\renewcommand{\ker}{\mathop{\rm Ker}\nolimits}
\DeclareMathOperator{\pgcd}{pgcd}

\newcommand{\reel}{\mathop{\mathrm{\reelname}}\nolimits}

\DeclareMathOperator{\Hom}{Hom}

\DeclareMathOperator{\DR}{DR}

\DeclareMathOperator{\SP}{SP}

\DeclareMathOperator{\diag}{diag}

\DeclareMathOperator{\image}{image}

\def\mod{{\mathrm{mod}}}

\newcommand{\loccit}{\emph{loc\ptbl cit}}

\def\cf{\textit{cf.}\kern.3em}
\def\eg{\textit{e.g.},\ }
\def\ie{\textit{i.e.},\ }
\def\resp{\textit{resp.}\kern.3em}

\newcommand{\T}{\S\kern .15em }
\newcommand{\ptbl}{.\kern .15em }

\def\to{\mathchoice{\longrightarrow}{\rightarrow}{\rightarrow}{\rightarrow}}
\def\mto{\mathchoice{\longmapsto}{\mapsto}{\mapsto}{\mapsto}}
\def\hto{\mathrel{\lhook\joinrel\to}}
\def\implique{\mathchoice{\Longrightarrow}{\Rightarrow}{\Rightarrow}{\Rightarrow}}
\let\Implique\implique

\let\iff\ssi

\def\cC{\mathcal{C}}
\def\cD{\mathcal{D}}
\def\cE{\mathcal{E}}

\def\cH{\mathcal{H}}

\def\cN{\mathcal{N}}
\def\cO{\mathcal{O}}

\def\cS{\mathcal{S}}

\DeclareMathAlphabet{\mathcalmaigre}{U}{eus}{m}{n}
\def\ccH{\mathcalmaigre{H}}

\let\cX\ccX

\def\omegag{\boldsymbol{\omega}}

\def\bR{\boldsymbol{R}}

\def\CC{\mathbb{C}}

\def\NN{\mathbb{N}}
\def\PP{\mathbb{P}}
\def\QQ{\mathbb{Q}}
\def\RR{\mathbb{R}}
\def\ZZ{\mathbb{Z}}

\def\gF{\mathfrak{F}}

\newcommand{\Afu}{\mathbb{A}^{\!1}}

\newcommand{\Afuh}{\wh{\mathbb{A}}^{\!1}}

\newcommand{\an}{\mathrm{an}}

\newcommand{\cbbullet}{{\raisebox{1pt}{$\bbullet$}}}
\newcommand{\gr}{\mathrm{gr}}

\def\sfrac#1#2{{#1}/{#2}}

\def\defin{\mathrel{:=}}

\newcommand{\un}{\mathbf{1}}
\def\numero{n${}^\circ$\kern2pt}
\makeatother

\begin{document}
\frontmatter
\dedicatory{Dedicated to Yuri Manin}
\begin{abstract}
We give an explicit description of the canonical Frobenius structure attached (by the results of the first part of this article) to the polynomial $f(u_0,\dots,u_n)=w_0u_0+\cdots+w_nu_n$ restricted to the torus $U=\{(u_0,\dots,u_n)\in\mathbb{C}^{n+1}\mid\prod_iu_i^{w_i}=1\}$, for any family of positive integers $w_0,\dots,w_n$ such that $\gcd(w_0,\dots,w_n)=1$.
\end{abstract}

\begin{altabstract}
Nous donnons une description explicite de la structure de Frobenius associ\'ee (par les r\'esultats de la premi\`ere partie de cet article) au polyn\^ome $f(u_0,\dots,u_n)=w_0u_0+\cdots+w_nu_n$ restreint au tore $U=\{(u_0,\dots,u_n)\in\mathbb{C}^{n+1}\mid\prod_iu_i^{w_i}=1\}$ pour toute famille de poids $w_0,\dots,w_n$ tels que $\pgcd(w_0,\dots,w_n)=1$.
\end{altabstract}


\altkeywords{Syst\`eme de Gauss-Manin, r\'eseau de Brieskorn, vari\'et\'e de Frobenius}


\maketitle
\mainmatter

\section{Introduction}

\subsection{}
This article explains a detailed example of the general result developed in the first part \cite{D-S02a}. We were motivated by \cite{Barannikov00}, where S.~Barannikov describes a Frobenius structure attached to the Laurent polynomial $f(u_0,\dots,u_n)=u_0+\cdots+u_n$ restricted to the torus $U=\{(u_0,\dots,u_n)\in\CC^{n+1}\mid\prod_iu_i=1\}$, and shows that it is isomorphic to the Frobenius structure attached to the quantum cohomology of the projective space $\PP^n(\CC)$ (as defined \eg in \cite{Manin96}).

We will freely use the notation introduced in the first part \cite{D-S02a}. A reference like ``\T I.3.c'' will mean \cite[\T3.c]{D-S02a}.

In the following, we fix an integer $n\geq 2$ and positive integers $w_0,\dots,w_n$ such that $\gcd(w_0,\dots,w_n)=1$. It will be convenient to assume that $w_0\leq\cdots\leq w_n$.
We put
\begin{equation}\label{eq:mu}
\mu\defin \sum_{i=0}^nw_i.
\end{equation}

We will analyze the Gauss-Manin system attached to the Laurent polynomial
$$f(u_1,\dots,u_n)=w_0u_0+w_1u_1+\cdots+w_nu_n$$
restricted to the subtorus $U\subset(\CC^*)^{n+1}$ defined by the equation $$u_0^{w_0}\cdots u_n^{w_n}=1.$$
The case $\mu=n+1$ (and all $w_i$ equal to~$1$) was considered in \cite{Barannikov00}. We will not need any explicit use of Hodge Theory, as all computation can be made ``by hand''. We will use the method of \T I.3.c to obtain information concerning the Frobenius structure on any germ of universal deformation space of~$f$. As we have seen in \cite{D-S02a}, we have to analyze with some details the structure of the Gauss-Manin system and the Brieskorn lattice of $f$.

\subsection{}\label{subsec:param}
Fix a $\ZZ$-basis of $\{\sum_iw_ix_i=0\}\subset\ZZ^{n+1}$. It defines a $(n+1)\times n$ matrix $M$. Denote by $m_0,\dots,m_n$ the lines of this matrix. We thus get a parametrization of $U$ by $(\CC^*)^n$ by putting $u_i=v^{m_i}$ for $i=0,\dots,n$ and $v=(v_1,\dots,v_n)$. The vectors $m_0,\dots,m_n$ are the vertices of a simplex $\Delta\subset\ZZ^n$, which is nothing but the Newton polyhedron of $f$ when expressed in the coordinates $v$. Notice that the determinant of the $n\times n$ matrix $(m_0,\dots,\wh{m_i},\dots,m_n)$ is $\pm w_i$.

\begin{lemme}\label{lem:Mtame}
The Laurent polynomial $f$ is convenient and nondegenerate with respect to its Newton polyhedron.
\end{lemme}

\begin{proof}
The nondegeneracy follows from the linear independence of any $n$ distinct vectors among $m_0,\dots,m_n$. Clearly, $0$ is contained in the interior of~$\Delta$.
\end{proof}

We know then that $f$ is M-tame (\cf \T I.4) and we may therefore apply the results of \T I.2 to $f$. An easy computation shows that $f$ has $\mu$ simple critical points, which are the $\zeta(1,\dots,1)$ with $\zeta^\mu=1$, and thus $\mu$ distinct critical values $\mu\zeta$. We hence have $\mu(f)=\mu$.

\subsection{}
Denote by $\cS_w$ the \emph{disjoint} union of the sets
\[
\{\ell \mu/w_i\mid\ell=0,\dots,w_i-1\}\subset \QQ.
\]
Hence $\#\cS_w=\mu$. Number the elements of $\cS_w$ from $0$ to $\mu-1$ in an \emph{increasing} way, with respect to the usual order on $\QQ$. We therefore have $\cS_w=\{s_w(0),\dots,s_w(\mu-1)\}$ with $s_w(k)\leq s_w(k+1)$. In particular, we have
$$s_w(0)=\cdots=s_w(n)=0, \quad s_w(n+1)=\frac{\mu}{\max_iw_i}<n+1.$$
Moreover, using the involution $\ell\mapsto w_i-\ell$ for $\ell\geq1$, one obtains, for $k\geq n+1$, the relation
\begin{equation}\label{eq:sym}
s_w(k)+s_w(\mu+n-k)=\mu.
\end{equation}
We consider the function $\sigma_w:\{0,\dots,\mu-1\}\to\QQ$ defined by
\begin{equation}\label{eq:spectre}
\sigma_w(k)=k-s_w(k).
\end{equation}
Hence $\sigma_w(k)=k$ for $k=0,\dots,n$. That $s_w(\cbbullet)$ is increasing is equivalent to
\begin{equation}\label{eq:alphaleq}
\forall k=0,\dots,\mu-1,\quad \sigma_w(k+1)\leq \sigma_w(k)+1,
\end{equation}
where we use the convention $\sigma_w(\mu)=\sigma_w(0)=0$. We will prove:

\begin{mtheoreme}\label{th:spectrew}
The polynomial $\prod_{k=0}^{\mu-1}(S+\sigma_w(k))$ is equal to the spectral polynomial $\SP_f(S)$ attached to $f$ (\cf \T I.2.e).
\end{mtheoreme}

For instance, if we take the Laurent polynomial $f(u_0,\dots,u_n)$ on the torus $\prod u_i=1$, \ie $w_0=\cdots=w_n=1$, we get $\SP_f(S)=\prod_{k=0}^{n}(S+k)$.

Notice that the symmetry property \eqref{eq:sym} is a little bit more precise than the symmetry of the spectrum (\cf \cite{Sabbah96b}), which would say that, for any $j\in\{0,\dots,n\}$, $$\#\{k\mid\sigma_w(k)=j\}=\#\{k\mid\sigma_w(k)=n-j\}.$$
Indeed, for $k\in\{n+1,\dots,\mu-1\}$, \eqref{eq:sym} means that $\sigma_w(k)+\sigma_w(\mu+n-k)=n$ and we clearly have $\sigma_w(k)+\sigma_w(n-k)=n$ for $k=0,\dots,n$.

\subsection{}
Consider now the two $\mu\times\mu$ matrices
\begin{equation}\label{eq:A0Ainfty}
A_\infty=\diag\big(\sigma_w(0),\dots,\sigma_w(\mu-1)\big),\quad
\arraycolsep=5pt
A_0=\mu\begin{pmatrix}
0&&&&1\\
1&0&&&0\\
0&1&0&&0\\
\vdots&\ddots&\ddots&\ddots&\vdots\\
0&\cdots&0&1&0
\end{pmatrix}.
\end{equation}
Notice that $A_0$ is semisimple with distinct eigenvalues $\mu\zeta$, where $\zeta$ is a $\mu$-th primitive root of $1$. In the canonical basis $(e_0,\dots,e_{\mu-1})$ of the space $\CC^\mu$ on which these matrices act, consider the nondegenerate bilinear form $g$ defined by
\[
g(e_k,e_\ell)=
\begin{cases}
1&
\begin{cases}\text{if }0\leq k\leq n\text{ and }k+\ell=n,\\
\text{or if }n+1\leq k\leq \mu-1\text{ and }k+\ell=\mu+n,
\end{cases}\\
0&\text{otherwise},
\end{cases}
\]
with respect to which $A_\infty$ satisfies $A_\infty+{}^t\!A_\infty=n\id$. The data $(A_0,A_\infty,g,e_0)$ define (\cf \cite[Main Theorem p\ptbl188]{Dubrovin96}, see also \cite[\T II.3]{Manin96} or \cite[Th\ptbl VII.4.2]{Sabbah00}) a unique germ of semisimple Frobenius manifold at the point $(\mu,\mu\zeta,\dots,\mu\zeta^{\mu-1})\in\CC^\mu$.

The main result of this article is then:

\begin{mtheoreme}\label{th:main}
The canonical Frobenius structure on any germ of a universal unfolding of the Laurent polynomial $f(u_0,\dots,u_n)=\sum_iw_iu_i$ on $U$, as defined in \cite{D-S02a}, is isomorphic to the germ of universal semisimple Frobenius structure with initial data $(A_0,A_\infty,g,e_0)$ at the point $(\mu,\mu\zeta,\dots,\mu\zeta^{\mu-1})\in\CC^\mu$.
\end{mtheoreme}

\begin{remarque*}
It would be interesting to give an explicit description of the Gromov-Witten potential attached to this Frobenius structure.
\end{remarque*}

\section{The rational numbers $\sigma_w(k)$}
Let us be now more precise on the definition of $s_w(k)$. Define inductively the sequence $(a(k),i(k))\in \NN^{n+1}\times\{0,\dots,n\}$ by
\begin{align*}
a(0)&=(0,\dots,0),&i(0)&=0\\
a(k+1)&=a(k)+\un_{i(k)},&i(k+1)&=\min\{i\mid a(k+1)_i/w_i=\textstyle\min_ja(k+1)_j/w_j\}.
\end{align*}
It is immediate that $\module{a(k)}\defin\sum_{i=0}^na(k)_i=k$ and that, for $k\leq n+1$, we have $a(k)_i=1$ if $i<k$ and $a(k)_i=0$ if $i\geq k$. In particular, $a(n+1)=(1,\dots,1)$.

\begin{lemme}\label{lem:ak}
The sequence $(a(k),i(k))$ satisfies the following properties:
\begin{enumerate}
\item\label{lem:ak1}
for all $k\in \NN$, $\dfrac{a(k)_{i(k)}}{w_{i(k)}}\leq\dfrac{a(k+1)_{i(k+1)}}{w_{i(k+1)}} \leq\dfrac{a(k)_{i(k)}+1}{w_{i(k)}}$,
\item\label{lem:ak2}
$a(\mu)=(w_0,\dots,w_n)$ and for all $k\in\{0,\dots,\mu-1\}$, we have $a(k)_{i(k)}\leq w_{i(k)}-1$,
\item\label{lem:ak3}
the map $\{0,\dots,\mu-1\}\to\coprod_{i=0}^n\{0,\dots,w_i-1\}$, defined by $k\mapsto [i(k),a(k)_{i(k)}]$ is bijective.
\item\label{lem:ak4}
For $\ell\in\NN$, we have $i(k+\ell\mu)=i(k)$ and $a(k+\ell\mu)_{i(k)}=\ell w_{i(k)}+a(k)_{i(k)}$.
\end{enumerate}
\end{lemme}

We will then put $s_w(k)\defin\mu a(k)_{i(k)}/w_{i(k)}$. We have $s_w(k+\ell\mu)=\ell\mu+s_w(k)$ for $\ell\in\NN$.

\begin{Proof}
\begin{enumerate}
\item
By induction on $k$. If $i(k+1)=i(k)$, the result is clear. Otherwise, we have $\sfrac{a(k+1)_{i(k+1)}}{w_{i(k+1)}}=\sfrac{a(k)_{i(k+1)}}{w_{i(k+1)}}$ and the first inequality follows from the definition of $i(k)$. Similarly, the second inequality is given by the definition of $i(k+1)$.
\item
Let us first remark the implication
\[
a(k)_j\leq w_j\ \forall j\text{ and } \{j\mid a(k)_j<w_j\}\neq\emptyset\Implique a(k+1)_j\leq w_j\ \forall j.
\]
[Indeed, from the assumption we have $a(k)_{i(k)}< w_{i(k)}$, hence $a(k+1)_{i(k)}=a(k)_{i(k)}+1\leq w_{i(k)}$. For $j\neq i(k)$, $a(k+1)_j=a(k)_j\leq w_j$.] Therefore, there exists $k_0$ such that $a(k_0)=(w_0,\dots,w_n)$. Then $k_0=\module{a(k_0)}=\mu$. Moreover, by what we have just seen, we have $a(k)_{i(k)}< w_{i(k)}$ for $k<\mu$.
\item
The map does exist, after \eqref{lem:ak2}, is clearly injective, therefore bijective as the two sets have the same number of elements.
\item
We have $a(\mu)=(w_0,\dots,w_n)$, so that $i(\mu)=0$, and we may apply the reasoning of \eqref{lem:ak2} for $k=\mu,\dots,2\mu-1$, etc.
\hfill\qed
\end{enumerate}
\skipqed
\end{Proof}

\begin{remarque}
In general, the numbers $s_w(k)$ are rational. These are integers (hence the spectrum of $f$ is integral) if and only if the following condition holds:
\begin{equation}\label{eq:conditionw}
\forall i,\quad w_i\mid \mu=w_0+\cdots+w_n.
\end{equation}

Consider the simplex $\Delta(w)$ in $\RR^n$ obtained as the intersection of the hyperplane $\ccH=\big\{\sum_{i=0}^nw_ix_i=0\big\}\subset\RR^{n+1}$ with the half spaces $x_i\geq-1$. Fix also the lattice $\ccH_\ZZ=H\cap\ZZ^{n+1}$. Then Condition \eqref{eq:conditionw} is equivalent to the condition that the vertices of $\Delta(w)$ are contained in the lattice $\ccH_\ZZ$. In other words, $\Delta(w)$ is a \emph{reflexive simplex} in the sense of Batyrev \cite{Batyrev94}. For instance, if $n=3$, one finds the following possibilities for $w_i$ (up to a permutation):

\noindent
\begin{small}
\mbox{}\hfill
\begin{minipage}[t]{4.8cm}
\[
\arraycolsep5pt
\def\arraystretch{1.2}
\begin{array}{|r|r|r|r||>{\bf}r|}
\hline
w_0&w_1&w_2&w_3&\mu\\
\hline
1&1&1&1&4\\
\hline
1&1&1&3&6\\
\hline
1&1&2&2&6\\
\hline
1&1&2&4&8\\
\hline
1&2&2&5&10\\
\hline
1&1&4&6&12\\
\hline
\end{array}
\]
\end{minipage}
\hfill
\begin{minipage}[t]{4.8cm}
\[
\arraycolsep5pt
\def\arraystretch{1.2}
\begin{array}{|r|r|r|r||>{\bf}r|}
\hline
w_0&w_1&w_2&w_3&\mu\\
\hline
1&2&3&6&12\\
\hline
1&3&4&4&12\\
\hline
1&2&6&9&18\\
\hline
1&4&5&10&20\\
\hline
1&3&8&12&24\\
\hline
2&3&10&15&30\\
\hline
1&6&14&21&42\\
\hline
\end{array}
\]
\end{minipage}
\hfill\mbox{}
\end{small}

\bigskip\noindent
For $n=4$, here are some examples (maybe not complete):
\par\nopagebreak
\begin{small}
\noindent
\mbox{}\hfill
\begin{minipage}[t]{4.8cm}
\[
\arraycolsep5pt
\def\arraystretch{1.2}
\begin{array}{|r|r|r|r|r||>{\bf}r|}
\hline
w_0&w_1&w_2&w_3&w_4&\mu\\
\hline
1&1&1&1&2&6\\
\hline
1&1&2&2&2&8\\
\hline
1&1&1&1&4&8\\
\hline
1&1&1&3&3&9\\
\hline
1&1&1&2&5&10\\
\hline
2&2&2&3&3&12\\
\hline
1&1&3&3&4&12\\
\hline
1&1&2&2&6&12\\
\hline
1&1&1&3&6&12\\
\hline
1&1&3&5&5&15\\
\hline
\end{array}
\]
\end{minipage}
\hfill
\begin{minipage}[t]{4.8cm}
\[
\arraycolsep5pt
\def\arraystretch{1.2}
\begin{array}{|r|r|r|r|r||>{\bf}r|}
\hline
w_0&w_1&w_2&w_3&w_4&\mu\\
\hline
1&1&2&4&8&16\\
\hline
1&1&4&4&10&20\\
\hline
1&1&4&6&12&24\\
\hline
1&1&2&8&12&24\\
\hline
1&1&3&10&15&30\\
\hline
1&1&4&12&18&36\\
\hline
1&1&8&10&20&40\\
\hline
1&1&6&16&24&48\\
\hline
1&1&8&20&30&60\\
\hline
1&2&12&15&30&60\\
\hline
\end{array}
\]
\end{minipage}
\hfill\mbox{}
\end{small}
\end{remarque}

\section{The Gauss-Manin system}
The Gauss-Manin system $G$ of the Laurent polynomial $f$ is a module over the ring $\CC[\tau,\tau^{-1}]$. It is defined as in \T I.2.c:
\[
G=\Omega^n(U)[\tau,\tau^{-1}]\big/(d-\tau df\wedge) \Omega^{n-1}(U)[\tau,\tau^{-1}].
\]
Put $\theta =\tau^{-1}$. The Briekorn lattice $G_{0}=\image(\Omega^{n}(U)[\theta]\to G)$
is a free $\CC[\theta]$-module of rank $\mu$ because, by Lemma~\ref{lem:Mtame}, $f$ is convenient and nondegenerate (\loccit.). We will consider the increasing filtration $G_p=\tau^pG_0$ ($p\in\ZZ$). Let $\omega_0$ be the $n$-form on $U$ defined by
\[
\omega_0=\dfrac{\frac{du_0}{u_0}\wedge\cdots\wedge\frac{du_n}{u_n}}{d\big(\prod_iu_i^{w_i}\big)}\Big|_{\prod_iu_i^{w_i}=1}.
\]
Let $v\mapsto u=v^m$ be a parametrization of $U$ as in \T\ref{subsec:param}. The form $\omega_0$ can be written as \hbox{$\omega_0=\pm\frac{dv_1}{v_1}\wedge\cdots\wedge\frac{dv_n}{v_n}$}. The Gauss-Manin system $G$ is then identified with the $\CC[\tau,\tau^{-1}]$-module (putting $v=(v_1,\dots,v_n)$)
\[
\CC[v,v^{-1},\tau,\tau^{-1}]\big/\{v_j\partial_{v_j}(\varphi_j)-\tau(v_j\partial_{v_j}f)\varphi_j\mid\varphi_j\in\CC[v,v^{-1},\tau,\tau^{-1}],\;j=1,\dots,n\}.
\]
It comes equipped with an action of $\partial_\tau$: if $\psi\in\CC[v,v^{-1}]$, let $[\psi]$ denote its class in $G$; then $\partial_\tau[\psi]=[-f\psi]$ (this does not depend on the representative of the class). Using the coordinate $\theta$, we have $\theta^2\partial_\theta[\psi]=[f\psi]$; this action is extended in the usual way to Laurent polynomials in $\tau$ with coefficients like $[\psi]$.

It is convenient to use the coordinates $u=(u_0,\dots,u_n)$. Then the previous quotient is written as
\[
\CC[u,u^{-1},\tau,\tau^{-1}]\big/\big(I_w+\CC[u,u^{-1},\tau,\tau^{-1}](g(u)-1)\big),
\]
where we have put $g(u)=\prod_iu_i^{w_i}$ and $I_w$ is the $\CC[\tau,\tau^{-1}]$-submodule of $\CC[u,u^{-1},\tau,\tau^{-1}]$ consisting of the expressions
\begin{equation}\label{eq:gmu}
\sum_{i=0}^nm_{ji}\Big(u_i\frac{\partial}{\partial u_i}-\tau w_iu_i\Big) \varphi_j, \text{ with }\varphi_j\in\CC[u,u^{-1},\tau,\tau^{-1}],\ (j=1,\dots,n).
\end{equation}

Consider the sequence $(a(k),i(k))$ of Lemma \ref{lem:ak}, and for each $k=0,\dots,\mu$, put
\[
\omega_k=u^{a(k)}\omega_0\in G_0,
\]
Notice that $\omega_\mu=\omega_0$ and, using \eqref{eq:0mi} below, that $f\omega_0=\mu\omega_1$.

\begin{proposition}\label{prop:bernstein}
The classes of $\omega_0,\omega_1,\dots,\omega_{\mu-1}$ form a $\CC[\theta]$-basis $\omegag$ of $G_0$. Moreover, they satisfy the equation
\[
-\frac1\mu(\tau\partial_\tau+\sigma_w(k))\omega_k=\tau\omega_{k+1}\quad(k=0,\dots,\mu-1),
\]
and we have Bernstein's relation in $G$:
\[
\prod_{k=0}^{\mu-1}\big[-\tfrac1\mu(\tau\partial_\tau-s_w(k))\big]\cdot\omega_0=\tau^{\mu}\omega_0.
\]
The $V$-order of $\omega_k$ is equal to $\sigma_w(k)$ and $\omegag$ induces a $\CC$-basis of $\oplus_\alpha\gr_\alpha^V(G_0/G_{-1})$.
\end{proposition}

From Theorem~I.4.5, Lemma~I.4.3(3), and the symmetry \eqref{eq:sym}, we get
\begin{equation}\label{eq:Sym}
\text{for }k=0,\dots,\mu-1, \quad 0\leq\sigma_w(k)\leq n \quad\text{and}\quad \begin{cases}
\sigma_w(k)=0\implique k=0,\\
\sigma_w(k)=n\implique k=n.
\end{cases}
\end{equation}
This implies that, for any $\alpha\in{}]0,n[$, the length of a maximal subsequence \hbox{$\alpha,\alpha+1,\dots,\alpha+\ell$} of $\sigma_w(\cbbullet)$ is $\leq n$, and even $\leq n-1$ if $\alpha$ is an integer. In other words:

\begin{corollaire}\label{cor:max}
The length of any maximal nonzero integral (resp. nonintegral) constant subsequence of $s_w(\cbbullet)$ is $\leq n-1$ (resp. $\leq n$).
\end{corollaire}

The proposition also gives a Birkhoff normal form for~$G_0$:
\[
\theta^2\partial_\theta\omegag=\omegag A_0+\theta \omegag A_\infty
\]
with $A_0,A_\infty$ as in \eqref{eq:A0Ainfty}. The matrix $A_0$ is nothing but the matrix of multiplication by $f$ on $G_0/\theta G_0$ in the basis induced by~$\omegag$. Its eigenvalues are the critical values of~$f$, as expected. In the case where $\mu=n+1$ (and all $w_i$ equal to~$1$), we find that $A_\infty=\diag(0,1,\dots,n)$ and $A_0$ is as in \eqref{eq:A0Ainfty} with size $\mu=n+1$.

\begin{proof}[Proof of Proposition \ref{prop:bernstein}]
It will be convenient to select some coordinate, say $u_0$. Multiplying \eqref{eq:gmu} (applied to $\varphi_1=\cdots=\varphi_n=\varphi$) on the left by the inverse matrix of the matrix formed by the columns of $m_1,\dots,m_n$, one finds that, for any $\varphi\in\CC[u,u^{-1},\tau,\tau^{-1}]$, we have in $G$
\begin{equation}\label{eq:0mi}
\forall i=1,\dots,n\quad (\tfrac1{w_i}u_i\partial_{u_i}-\tfrac1{w_0}u_0\partial_{u_0})\varphi\omega_0 =\tau(u_i-u_0)\varphi\omega_0.
\end{equation}
Applying this to any monomial $\varphi=u^a$ and summing these equalities, we get the following relation for $j=0$, hence for any $j=0,\dots,n$ by a similar argument:
\begin{equation}\label{eq:Lj}
-\frac1\mu\big(\tau\partial_\tau+L_j(a)\big)u^a\omega_0=\tau u^{a+\un_j}\omega_0,
\end{equation}
where we put $L_j(a)=\sum_{i=0}^na_i-\mu a_j/w_j$. This is nothing but (I.4.12) in the present situation. Apply this for $a=a(k)$ and $j=i(k)$ ($k=0,\dots,\mu-1$) to get the first relation in the lemma (remark that $L_{i(k)}(a(k))=\sigma_w(k)$). Bernstein's relation for $\omega_0$ is then clear. Remark also that $\omega_k$ is given by
\[
\omega_k=\tau^{-k}\prod_{j=0}^{k-1}\Big[-\frac1\mu(\tau\partial_\tau-s_w(j))\Big]\cdot \omega_0.
\]
It is not difficult to derive from Bernstein's relation for $\omega_0$ a Bernstein relation for each $\omega_k$ and conclude that $\omega_k$ has $V$-order $\leq \sigma_w(k)$. [Notice also that, as $\sigma_w(k)=L_{i(k)}(a(k))=\max_jL_j(a(k))$, the order of $\omega_k$ with respect to the Newton filtration is $\leq \sigma_w(k)$; this is compatible with Theorem~I.4.5.]

Let us now show that $\omega_0,\dots,\omega_{\mu-1}$ generate $G_0$ as a $\CC[\theta]$-module. Notice that Bernstein's relation for $\omega_0$ implies that $\partial_\tau^\mu\omega_0\in\CC[\theta]\langle\omega_0,\dots,\partial_\tau^{\mu-1}\omega_0\rangle=\CC[\theta]\langle\omega_0,\dots,\omega_{\mu-1}\rangle$, and this also holds for $\partial_\tau^\ell\omega_0$ for $\ell\geq\mu$. It is therefore enough to show that $(f^\ell\omega_0)_{\ell\geq0}$ generate $G_0$ over $\CC[\theta]$. Write \eqref{eq:0mi} as
\begin{equation}\label{eq:ua}
u^{a+\un_i}\omega_0=\Big[u^au_0+\Big(\frac{a_i}{w_i}-\frac{a_0}{w_0}\Big)\theta u^a\Big]\omega_0.
\end{equation}
The Brieskorn lattice $G_0$ is generated over $\CC[\theta]$ by the $u_0^\ell\omega_0$ with $\ell\in\NN$: indeed, it is generated by the $u^a\omega_0$; then,
\begin{itemize}
\item
if $a_i\geq1$ for some $i\geq1$, one decreases $a_i$ to $0$ with \eqref{eq:ua};
\item
if $a_i\leq-1$ for some $i\geq1$, one iterates \eqref{eq:ua} $w_0$ times and use the relation $u^w\omega_0=\omega_0$ to express $u^b\omega_0$ (any $b$) as a sum (with constant coefficients) of terms $\theta^ku_0^\ell u^{b+w'}\omega_0$ and of $u^{b+w'+\un_i}$, with $k,\ell\geq0$ and $w'=(0,w_1,\dots,w_n)$; hence if $b_i<0$, there exists $r$ such that $b_i+rw_i\geq0$ and one iterates $r$ times the previous process to write $u^b\omega_0$ with terms $\theta^ku^a$, with $a_i\geq 1$, to reduce to the previous case;
\item
notice that, in both previous processes, we never decrease the degree in $u_0$; now, we are reduced to considering $u_0^\ell\omega_0$ with $\ell<0$; use once more the relation $u^{k w}u^b\omega_0=u^b\omega_0$ (for any $k\geq0$, any $b$) to replace $u_0^\ell$ with $u^a$ with $a_0,\dots, a_n\geq0$ and apply the first case.
\end{itemize}
A similar argument gives the result for the family $(f^\ell\omega_0)_{\ell\geq0}$. As $G_0$ is $\CC[\theta]$-free (\cf Remark~I.4.8 and \T I.2.c), we conclude that $\omegag$ is a $\CC[\theta]$-basis of $G_0$. [Instead of using Remark~I.4.8, one can directly conclude here that $G_0$ is $\CC[\theta]$-free of rank $\mu$ by showing first that $\omegag$ generates $G$ as a $\CC[\tau,\tau^{-1}]$-module.]

Remark also that $(\omega_0,\dots,f^{\mu-1}\omega_0)$ is another basis, but the differential equation does not take Birkhoff normal form in such a basis.

\medskip
We will now determine the $V$-filtration. Put $\omega'_k=\tau^{[\sigma_w(k)]}\omega_k$. Then $\omegag'$ is another $\CC[\tau,\tau^{-1}]$-basis of $G$. The $V$-order of $\omega'_k$ is $\leq \sigma_w(k)-[\sigma_w(k)]<1$. For $\alpha\in[0,1[$, put
\begin{align*}
U_\alpha G&=\CC[\tau]\langle\omega'_k\mid \sigma_w(k)-[\sigma_w(k)]\leq\alpha\rangle+\tau\CC[\tau]\langle\omega'_k\mid \sigma_w(k)-[\sigma_w(k)]>\alpha\rangle\\
U_{<\alpha}G&=\CC[\tau]\langle\omega'_k\mid \sigma_w(k)-[\sigma_w(k)]<\alpha\rangle+\tau\CC[\tau]\langle\omega'_k\mid \sigma_w(k)-[\sigma_w(k)]\geq\alpha\rangle,
\end{align*}
and $U_{\alpha+p}G=\tau^pU_\alpha G$ (\resp $U_{<\alpha+p}G=\tau^pU_{<\alpha} G$) for any $p\in\ZZ$. We then have
\[
U_\alpha G=\CC\langle\omega'_k\mid \sigma_w(k)-[\sigma_w(k)]=\alpha\rangle+ U_{<\alpha}G.
\]
Notice that, according to the formula for $\omega_k$, the elements $\omega'_k$ satisfy
\begin{equation}\label{eq:tauomega'}
\begin{split}
-\frac1\mu\big(\tau\partial_\tau+\sigma_w(k)-[\sigma_w(k)]\big)\omega'_k&= \tau^{[\sigma_w(k)]+1-[\sigma_w(k+1)]}\omega'_{k+1}\\
&=\tau^{\lceil s_w(k+1)\rceil-\lceil s_w(k)\rceil}\omega'_{k+1},
\end{split}
\end{equation}
with $\lceil s\rceil\defin -[-s]$. Recall that the sequence $(s_w(k))$, hence the sequence $(\lceil s_w(k)\rceil)$, is increasing. If $\lceil s_w(k+1)\rceil>\lceil s_w(k)\rceil$, then
\[
\big(\tau\partial_\tau+\sigma_w(k)-[\sigma_w(k)]\big)\omega'_k\in U_{<0}G.
\]
Otherwise, we have
\[
\lceil s_w(k)\rceil-s_w(k)\geq\lceil s_w(k+1)\rceil-s_w(k+1),
\]
\ie
\[
\sigma_w(k)-[\sigma_w(k)]\geq \sigma_w(k+1)-[\sigma_w(k+1)],
\]
and we conclude that $U_\alpha G$ is stable under $\tau\partial_\tau$ and that $\tau\partial_\tau+\alpha$ is nilpotent on $\gr_\alpha^UG$. The filtration $U_\bbullet G$ satisfies then the characterizing properties of $V_\bbullet G$, hence is equal to it.

\medskip
We may now compute $G_p\cap V_\alpha$ for $p\in\ZZ$ and $\alpha\in[0,1[$. Any element of $G_p\cap V_\alpha$ decomposes uniquely as $\sum_{k=0}^{\mu-1}p_k(\tau)\omega'_k$, with
\[
p_k(\tau)\in\begin{cases}
\tau^{p-[\sigma_w(k)]}\CC[\tau^{-1}]\cap\CC[\tau]&\text{if } \sigma_w(k)-[\sigma_w(k)]\leq\alpha,\\
\tau^{p-[\sigma_w(k)]}\CC[\tau^{-1}]\cap\tau\CC[\tau]&\text{if } \sigma_w(k)-[\sigma_w(k)]>\alpha
\end{cases}
\]
It follows that
\begin{equation}\label{eq:GpVa}
G_p\cap V_\alpha= \sum_{k\mid\sigma_w(k)=\alpha+p}\hspace*{-5mm}\CC\cdot\omega'_k+G_p\cap V_{<\alpha} +G_{p-1}\cap V_\alpha,
\end{equation}
and therefore $\gr_p^G\gr_\alpha^VG$ is generated by the classes of $\omega'_k$ with $\sigma_w(k)=\alpha+p$. These classes form a basis of $\gr_p^G\gr_\alpha^VG$, as $\dim\oplus_p\oplus_{\alpha\in[0,1[}\gr_p^G\gr_\alpha^VG=\mu$. This gives the last statement of the proposition.
\end{proof}

For $\alpha\in [0,1[$, let $\omega'_k$ be such that $\lceil s_w(k)\rceil-s_w(k)=\alpha$ and denote by $[\omega'_k]$ the class of $\omega'_k$ in $H_\alpha\defin\gr_\alpha^VG$. After \eqref{eq:alphaleq} we have:
\begin{equation}\label{eq:omegaprime}
-\frac1\mu(\tau\partial_\tau+\alpha)[\omega'_k]=
\begin{cases}
0&\text{if }s_w(k+1)>s_w(k),\\
[\omega'_{k+1}]&\text{if }s_w(k+1)=s_w(k).
\end{cases}
\end{equation}
It follows that the \emph{primitive elements} relative to the nilpotent operator induced by $(-1/\mu)(\tau\partial_\tau+\alpha)$ on $H_\alpha$ are
the elements $[\omega'_k]$ such that
\[
k\geq n+1,\quad \lceil s_w(k)\rceil-s_w(k)=\alpha\quad\text{and}\quad s_w(k-1)<s_w(k)
\]
and, if moreover $\alpha=0$, the element $[\omega'_0]=[\omega_0]$.

Therefore, the Jordan blocks of $(-1/\mu)(\tau\partial_\tau+\alpha)$ on $H_\alpha$ are in one-to-one correspondence with the maximal constant sequences in $\cS_w$, and the corresponding sizes are the same. All Jordan blocks, except that of $[\omega_0]$ if $\alpha=0$, have thus size $\leq n$, and even $\leq n-1$ if $\alpha$ is an integer (\cf Cor\ptbl\ref{cor:max}). Recall also (\cf \cite{Sabbah96a,Sabbah96b}) that $H\defin\oplus_{\alpha\in[0,1[}H_\alpha$ may be identified with the relative cohomology space $H^n(U,f^{-1}(t))$ for $\module{t}\gg0$, that $H_\alpha$ corresponds to the generalized eigenspace of the monodromy corresponding to the eigenvalue $\exp2i\pi\alpha$, and that the unipotent part of the monodromy operator $T$ is equal to $\exp2i\pi N$ with $N\defin-(\tau\partial_\tau+\alpha)$.

\begin{exemple}
Take $n=4$ and $w_0=1$, $w_1=2$, $w_2=12$, $w_3=15$ and $w_4=30$, so that $\mu=60$. Then the only possible $\alpha$ is $0$ and $N$ has one Jordan block of size~$5$, $3$~blocks of size~$3$, $13$~blocks of size~$2$ and $20$~blocks of size~$1$. On the other hand, if $\mu=n+1$ (and all $w_i$ equal to~$1$), the only possible $\alpha$ is $0$ and $N$ has only one Jordan block (of size $n+1$).
\end{exemple}

\section{Poincar\'e duality and higher residue pairings}
Consider on $\CC[\tau,\tau^{-1}]$ the ring involution induced by $\tau\mapsto-\tau$. We will set $\ov{p(\tau)}\defin p(-\tau)$ (there is no complex conjugation involved here). Given a $\CC[\tau,\tau^{-1}]$-module $G$, we denote by $\ov G$ the $\CC$-vector space $G$ equipped with the new module structure $p(\tau)\cdot g=p(-\tau)g$. For convenience, we denote by $\ov g$ the elements of $\ov G$. The $\CC[\tau,\tau^{-1}]$-structure of $\ov G$ is therefore given by the rule: $\ov{p(\tau)}\ov{g}=\ov{p(\tau)g}$.

If $G$ is moreover equipped with a connection, \ie with a compatible action of $\partial_\tau$, then so is $\ov G$ and we have $\partial_\tau\ov{g}\defin\ov{-\partial_\tau g}$. Notice that $\tau\partial_\tau\ov g=\ov{\tau\partial_\tau g}$.

Duality for $\cD$-modules gives (\cf \cite{Sabbah96b}) the existence of a nondegenerate $\CC[\tau,\tau^{-1}]$-bilinear pairing
\[
S:G\ootimes_{\CC[\tau,\tau^{-1}]}\ov G\longrightarrow \CC[\tau,\tau^{-1}]
\]
satisfying the following properties:
\begin{enumerate}
\item\label{SLeibniz}
$\dfrac{dS(g',\ov{g''})}{d\tau}=S(\partial_\tau g',\ov{g''})+S(g',\partial_\tau\ov{g''})=S(\partial_\tau g',\ov{g''})-S(g',\ov{\partial_\tau g''})$,
\par\noindent
(equivalently, $\tau\partial_\tau S(g',\ov{g''})=S(\tau\partial_\tau g',\ov{g''})+S(g',\ov{\tau\partial_\tau g''})$),
\item\label{SV}
$S$ sends $V_0\otimes\ov{V_{<1}}$ in $\CC[\tau]$,
\item\label{SG}
$S$ sends $G_0\otimes \ov{G_0}$ in $\theta^n\CC[\theta]=\tau^{-n}\CC[\tau^{-1}]$,
\item\label{Ssym}
$S(g'',\ov{g'})=(-1)^n\ov{S(g',\ov{g''})}$ (this reflects the $(-1)^n$-symmetry of the Poincar\'e duality on $U$).
\end{enumerate}

Notice that \eqref{SLeibniz} means that $S$ is a horizontal section of the $\CC[\tau,\tau^{-1}]$-module $\Hom_{\CC[\tau,\tau^{-1}]}(G\otimes\ov G,\CC[\tau,\tau^{-1}])$ equipped with its natural connection, or also that $S$ is a $\CC[\tau]\langle\partial_\tau\rangle$-linear morphism $\ov G\to G^*$, if one endows $G^*=\Hom_{\CC[\tau,\tau^{-1}]}(G,\CC[\tau,\tau^{-1}])$ with its natural connection. Therefore, \eqref{SV} follows from \eqref{SLeibniz} because any $\CC[\tau]\langle\partial_\tau\rangle$-linear morphism is strict with respect to the Malgrange-Kashiwara filtrations $V$ and we have
\[
V_\beta(G^*)=\Hom_{\CC[\tau]}(V_{<-\beta+1}G,\CC[\tau])
\]
(\cf \cite{Sabbah96b}).

In the case of singularities, this corresponds to the ``higher residue pairings'' of K.~Saito \cite{KSaito83}. The link with Poincar\'e duality is explained in \cite{MSaito89}.

For our Laurent polynomial $f$, we will recover in an elementary way the existence of such a pairing $S$ satisfying the previous properties. More precisely, we have:

\begin{lemme}\label{lem:dualite}
There exists a unique (up to a nonzero constant) nondegenerate pairing~$S$ satisfying Properties \eqref{SLeibniz}, \eqref{SV}, \eqref{SG}. It is given by the formula:
\[
S(\omega_k,\ov{\omega_\ell})=
\begin{cases}
S(\omega_0,\ov{\omega_n})&
\begin{cases}\text{if }0\leq k\leq n\text{ and }k+\ell=n,\\
\text{or if }n+1\leq k\leq \mu-1\text{ and }k+\ell=\mu+n,
\end{cases}\\
0&\text{otherwise}.
\end{cases}
\]
Moreover, for any $k,\ell$, $S(\omega_k,\ov{\omega_\ell})$ belongs to $\CC\tau^{-n}$ and $S$ satisfies \eqref{Ssym}.
\end{lemme}

\begin{proof}
Assume that a pairing $S$ satisfying \eqref{SLeibniz}, \eqref{SV}, \eqref{SG} exists. For $k,\ell=0,\dots,\mu-1$, we have $S(\omega_k,\ov{\omega_\ell})\in\tau^{-n}\CC[\tau^{-1}]$ by \eqref{SG} and $S(\omega_0,\ov{\omega_\ell})\in\tau^{-[\sigma_w(\ell)]}\CC[\tau]$ by \eqref{SV}. Therefore, $S(\omega_0,\ov{\omega_\ell})\neq0$ implies $[\sigma_w(\ell)]\geq n$, and if $[\sigma_w(\ell)]=n$, we have $S(\omega_0,\ov{\omega_\ell})\in\CC\tau^{-n}$. But we know by \eqref{eq:Sym} that
$$[\sigma_w(\ell)]
\begin{cases}
{}<n &\text{if }\ell\neq n,\\
{}=n&\text{if }\ell=n.
\end{cases}
$$
Therefore, $S(\omega_0,\ov{\omega_\ell})=0$ if $\ell\neq n$ and $S(\omega_0,\ov{\omega_n})\in\CC\tau^{-n}$.

Notice also that we have by \eqref{SLeibniz} and Proposition \ref{prop:bernstein}:
\begin{multline}\label{eq:induc}
-\frac1\mu(\tau\partial_\tau+n)S(\omega_k,\ov{\omega_\ell})\\
=\tau\big[S(\omega_{k+1},\ov{\omega_\ell})- S(\omega_k,\ov{\omega_{\ell+1}}) \big]+\frac{\sigma_w(k)+\sigma_w(\ell)-n}{\mu}S(\omega_k,\ov{\omega_\ell}),
\end{multline}
if we put as above $\omega_\mu=\omega_0$.

Argue now by induction for $k<n$: as $S(\omega_k,\ov{\omega_\ell})\in\CC\tau^{-n}$, the LHS in \eqref{eq:induc} vanishes. This shows that $S(\omega_{k+1},\ov{\omega_\ell})=0$ if $\ell\neq n-k,n-1-k$. Moreover, if $\ell=n-k$, we have \hbox{$\sigma_w(k)+\sigma_w(\ell)-n=0$}, hence $S(\omega_{k+1},\ov{\omega_{n-k}})=0$. Last, we have $S(\omega_{k+1},\ov{\omega_{n-1-k}})=S(\omega_{k},\ov{\omega_{n-k}})$.

Argue similarly for $k\geq n+1$.
\end{proof}

Notice that, if $A_\infty^*$ denotes the adjoint of $A_\infty$ with respect to $S$, then $A_\infty+A_\infty^*=n\id$, \ie $A_\infty-(n/2)\id$ is skewsymmetric with respect to $S$.

\section{M.~Saito's solution to the Birkhoff problem}
One step in constructing the Frobenius structure associated to $f$ consists in solving Birkhoff's problem for the Brieskorn lattice $G_0$ in the Gauss-Manin system $G$, that is, in finding a $\CC[\tau]$-lattice $E$ of $G$, which glues with $G_0$ to a trivial vector bundle on $\PP^1$. Recall (\cf \cite[App\ptbl B]{D-S02a} for what follows) that there is a one-to-one correspondence between such lattices $E$ which are \emph{logarithmic}, and decreasing filtrations $\oplus_{\alpha\in[0,1[}H_\alpha^\cbbullet$ of $H=\oplus_{\alpha\in[0,1[}H_\alpha$ which are stable under $N$ and which are \emph{opposite} to the filtration
$$
G_p(H)\defin \ooplus_{\alpha\in[0,1[}(G_p\cap V_\alpha)/(G_p\cap V_{<\alpha})=\CC\langle[\omega'_k]\mid[\sigma_w(k)]\leq p\rangle\quad\text{after \eqref{eq:GpVa}}.
$$
This is analogous to \cite[Th\ptbl3.6]{MSaito89}.

In \cite[Lemma 2.8]{MSaito89}, M.~Saito defines a canonical decreasing filtration $H^\cbbullet_{\textup{Saito}}$ in terms of the monodromy filtration $M_\bbullet$ of the nilpotent endomorphism $2i\pi N$ of $H$ and of the filtration conjugate to $G_p(H)$, the conjugation being taken with respect to the real structure on $H$ coming from the identification with $H^n(U,f^{-1}(t))$. This defines therefore a canonical solution to Birkhoff's problem for $G_0$.

Consider now the decreasing filtration $H^\cbbullet$ of $H$ explicitly defined by
\begin{equation}\label{eq:Hp}
H^p=\CC\langle[\omega'_k]\mid [\sigma_w(k)]\geq p\rangle,
\end{equation}
where $[\,]$ denotes the class in $H=\oplus_{\alpha\in[0,1[}V_\alpha G/V_{<\alpha}G$. Then $H^\cbbullet$ is opposite to $G_\bbullet(H)$. It satisfies
\[
H^0=H,\ H^{n+1}=0,\quad NH^p\subset H^{p+1}
\]
and, for $k=0,\dots,\mu-1$ and $\alpha\in[0,1[$,
\[
(H_\alpha^p)^\perp=
\begin{cases}
H_{1-\alpha}^{n-p}&\text{if }\alpha\neq0,\\
H_0^{n+1-p}&\text{if }\alpha=0,
\end{cases}
\]
where ${}^\perp$ means taking the orthogonal with respect to the symmetric bilinear form $g$ on $H$ induced by $S$. If $\mu=n+1$ (and all $w_i$ equal to~$1$), then $H^p=M_{n-2p}$ (this implies that the mixed Hodge structure on $H$ is ``Hodge-Tate'').

\begin{proposition}\label{prop:Hodgeoppose}
The filtration $H^\cbbullet$ is equal to the opposite filtration $H^\cbbullet_{\textup{Saito}}$. The associated logarithmic lattice is $E\defin\CC[\tau]\langle\omega_0,\dots,\omega_{\mu-1}\rangle$.
\end{proposition}

\begin{proof}
Let us begin with the second statement. The lattice $E$ introduced in the proposition is logarithmic, by Proposition~\ref{prop:bernstein}. A computation analogous to that of $G_p\cap V_\alpha$ shows that the filtration $\tau^p E\cap V_\alpha/V_{<\alpha}$ of $H_\alpha$ is equal to $H^p_\alpha$. Therefore, $E$ is the logarithmic lattice corresponding to $H^\cbbullet$ by the correspondence recalled above.

For the first statement, put $F^\cbbullet (H)=G_{n-\bbullet}(H)$. This is a decreasing filtration. Consider also the increasing filtration
$$
W_\bbullet(H_\alpha)=\begin{cases}
M_{\bbullet-n-1}(H_\alpha)&\text{if }\alpha\neq0,\\
M_{\bbullet-n}(H_0)&\text{if }\alpha=0,
\end{cases}
$$
where $M_\bbullet(H)$ denotes the monodromy filtration of the nilpotent endomorphism $2i\pi N$ on~$H$. Recall that $W_\bbullet(H)$ is defined over $\RR$ (even over $\QQ$) as $2i\pi N$ is so. Then the opposite filtration given by M.~Saito is
\[
H^\bbullet_{\textup{Saito}}=\sum_q\ov F^q\cap W_{n+q-\bbullet}(H),
\]
where $\ov E$ denotes the conjugate of the subspace $E$ of $H$ with respect to the complex conjugation coming from the identification
$$H\isom H^n(U,f^{-1}(t),\CC)=\CC\otimes_\RR H^n(U,f^{-1}(t),\RR).$$
We therefore need to give a description of the conjugation in term of the basis $[\omega'_k]$.

Let $k_0$ be such that $[\omega'_{k_0}]$ is a primitive element with respect to $N$, and denote by $\nu_{k_0}$ its weight. Then $N^{\nu_{k_0}+1}[\omega'_{k_0}]=0$. For $j=0,\dots,\nu_{k_0}$, put $k=k_0+j$. Then $[\omega'_k]=(\tfrac1\mu N)^j[\omega'_{k_0}]$ has order $\nu_{k_0}+1-j$ with respect to $N$, and weight $\nu_k\defin\nu_{k_0}-2j$. Moreover, we have $\sigma_w(k)=\sigma_w(k_0)+j$, as $j\mto s_w(k_0+j)$ is constant. The space $B_{k_0}\defin\langle N^j[\omega'_{k_0}]\mid j=0,\dots,\nu_{k_0}\rangle$ is a Jordan block of $N$.

Assume that $k_0\geq n+1$. Then $[\omega'_{\mu+n-k_0}]$ is primitive with respect to ${}^t\!N$, hence $[\omega'_{\mu+n-k_0-\nu_{k_0}}]$ is primitive with respect to $N$. It will be convenient to put $\ov k_0=\mu+n-k_0-\nu_{k_0}$ and, for $k=k_0+j$ with $j=0,\dots,\nu_{k_0}$, $\ov k=\ov k_0+j$. We therefore have $\ov k=\mu+n-k-\nu_k$. Notice that, for such a $k$, we have $s_w(k)=s_w(k_0)=\mu-s_w(\ov k_0)=\mu-s_w(\ov k)$.
We also have $\sigma_w(\mu+n-k)-\nu_k=\sigma_w(\ov k)$ if $k\geq n+1$.

For $k\in[0,n]$, we simply put $\ov k=k$.

The proof of the following lemma will be given in \T\ref{sec:topology}.

\begin{lemme}\label{lem:conjug}
For $k_0\geq n+1$, the conjugate of the Jordan block $B_{k_0}$ is the Jordan block $B_{\ov k_0}$, and $B_0$ is self-conjugate.
\end{lemme}

It follows from this lemma that, for $k$ as above, we have
\begin{equation}\label{eq:conjomega}
\ov{[\omega'_k]}=\sum_{\ell=k}^{k_0+\nu_0}a_\ell[\omega'_{\ov\ell}]
\end{equation}
with $a_k\neq0$.

Let us now end the proof of Proposition \ref{prop:Hodgeoppose}. We have
\[
F^q\cap W_{n+q-p}=G_{n-q}\cap M_{q-p(-1)} =\langle[\omega'_k]\mid[\sigma_w(k)]\leq n-q \text{ and }\nu_k\leq q-p(-1)\rangle,
\]
where $(-1)$ is added if $\sigma_w(k)\not\in\ZZ$ and not added otherwise. Therefore,
\[
\sum_qF^q\cap W_{n+q-p}=\langle[\omega'_k]\mid[\sigma_w(k)]+\nu_k\leq n-p(-1)\rangle.
\]
Remark now that, if $k\geq n+1$,
\begin{align*}
[\sigma_w(k)]+\nu_k\leq n-p(-1)&\iff [n-\sigma_w(\mu+n-k)]+\nu_k\leq n-p(-1)\\
&\iff\lceil\sigma_w(\mu+n-k-\nu_k)\rceil\geq p(+1)\\
&\iff[\sigma_w(\mu+n-k-\nu_k)]\geq p\\
&\iff [\sigma_w(\ov k)]\geq p.
\end{align*}
Arguing similarly for $k\leq n$, we conclude from Lemma \ref{lem:conjug} and \eqref{eq:Hp} that
\begin{align}
\smash{\sum_q}\ov F^q\cap W_{n+q-p}&=\langle\ov{[\omega'_k]}\mid[\sigma_w(k)]+\nu_k\leq n-p(-1)\rangle\notag\\
&=\langle[\omega'_{\ov k}]\mid[\sigma_w(\ov k)]\geq p\rangle\label{eq:conjB}\\
&=H^p.\notag
\end{align}
Notice that \eqref{eq:conjB} follows from \eqref{eq:conjomega}, as $\sigma_w$ is increasing on each $B_{\ov k_0}$.
\end{proof}

\section{Some topology of $f$ and proof of Lemma \ref{lem:conjug}} \label{sec:topology}

\subsection{Lefschetz thimbles}
Denote by $\Delta$ the subset $(\RR_+^*)^{n+1}\cap U\subset U$, defined by $u_i>0$ for $i=0,\dots,n$. The restriction $f_{|\Delta}$ of $f$ to $\Delta$ takes values in $[\mu,+\infty[$ and has only one critical point (which is a Morse critical point of index $0$), namely $(1,\dots,1)$, with critical value equal to $\mu$. Notice also that $f_{|\Delta}$ is proper. Therefore, $\Delta$ is a Lefschetz thimble for $f$ with respect to the critical point $(1,\dots,1)$. Other Lefschetz thimbles at $\zeta^\ell(1,\dots,1)$ are $\zeta^\ell\Delta$ (\hbox{$\ell=0,\dots,\mu-1$}).

Fix $\tau\neq0$. The morphisms
\[
H_n(U,\reel\tau f>C';\QQ)\to H_n(U,\reel\tau f>C;\QQ)
\]
for $C'>C$ are isomorphisms if $C$ is big enough. We denote by \hbox{$H_n(U,\reel\tau f\gg0)$} the limit of this inverse system. This is the germ at $\tau$ of a local system $\cH$ of rank $\mu$ on $\CC^*=\{\tau\neq0\}$.
Notice that $\Delta$ defines a nonzero element of the germ $\cH_\tau$ at any $\tau$ with $\reel\tau>0$, \ie a section $\Delta(\tau)$ of $\cH$ on $\{\reel\tau>0\}$. Therefore, it defines in a unique way a multivalued section of $\cH$ on~$\CC^*$.

Let $\epsilon>0$ be small enough. As $f$ is a $C^\infty$ fibration over the open set \hbox{$\CC\moins\{\mu\zeta^\ell\mid\ell=0,\dots,\mu-1\}$}, it is possible to find a basis of sections \hbox{$\Delta_0(\tau),\dots,\Delta_{\mu-1}(\tau)$} of $\cH$ on the open set
$$
\cS=\{\tau=\module{\tau}e^{2i\pi\theta}\mid\theta\in{}]\epsilon-1,\epsilon[\}
$$
in such a way that, for any $\ell\in\{0,\dots,\mu-1\}$ and $\module{\tau}>0$, we have $\Delta_\ell(\zeta^{-\ell}\module{\tau})=\zeta^\ell\Delta$. Of course, this basis extends as a basis of multivalued sections of $\cH$ on~$\CC^*$.

\begin{figure}[htb]
\begin{tabular}{ccc}
\includegraphics[scale=.65]{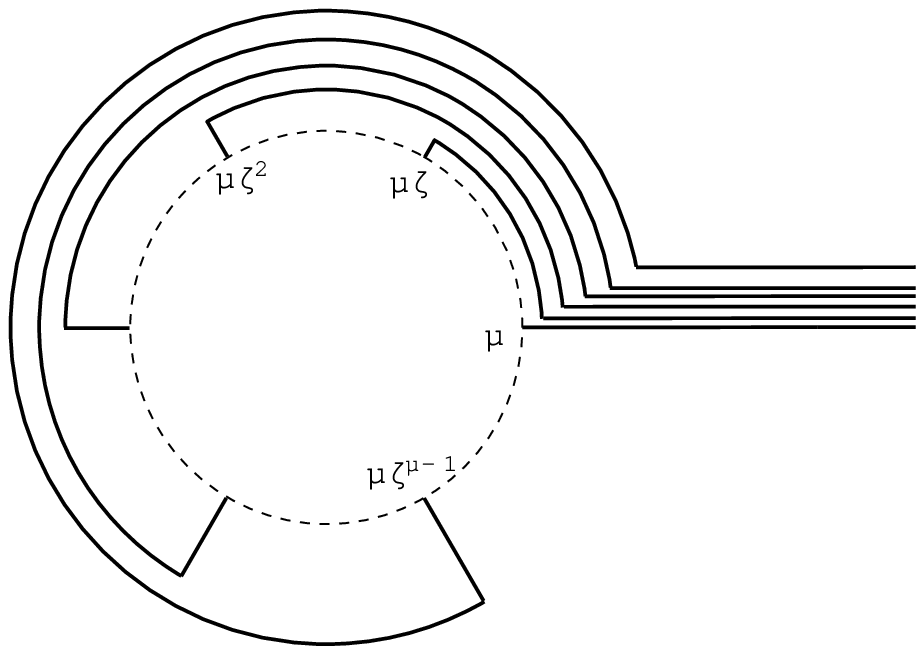}&\hspace*{1cm}&\includegraphics[scale=.65]{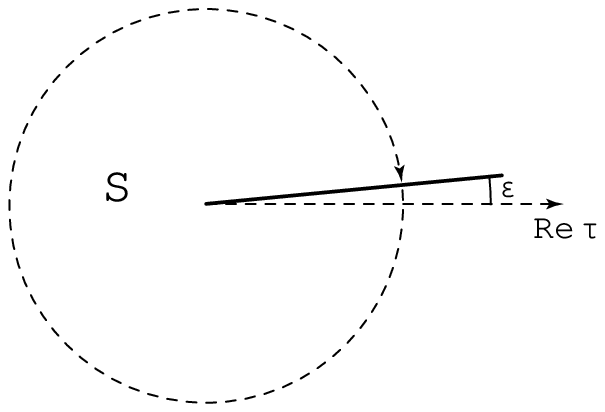}\\
The images $f(\Delta_\ell(\module{\tau}))$&&The domain $\cS$
\end{tabular}
\end{figure}

\subsection{Integrals along Lefschetz thimbles}
Let $\eta\in G$ and let $\wt\eta$ be a representative of $\eta$ in $\Omega^n(U)[\tau,\tau^{-1}]$. Then the function
\[
\cS\ni\tau\mto\int_{\Delta_\ell(\tau)}e^{-\tau f}\wt\eta
\]
only depends on $\eta$ and is holomorphic on $\cS$. It is denoted by $\varphi_{\ell,\eta}(\tau)$. Moreover, we have
\[
\frac{\partial\varphi_{\ell,\eta}(\tau)}{\partial\tau}=\varphi_{\ell,\partial_\tau\eta}(\tau).
\]
\let\oldlog\log\def\log{\tfrac{1}{2i\pi}\oldlog}
It follows that, if $\eta\in V_\alpha G$, then
\begin{equation}\label{eq:psieta}
\varphi_{\ell,\eta}(\tau)=\tau^{-\alpha}\Big[\sum_{m=0}^{m_\eta}c_{\eta,\ell}^{(m)}\frac{(-\log\tau)^m}{m!}+o(1)\Big]\defin\tau^{-\alpha}\big[\psi_{\eta,\ell}(\tau)+o(1)\big],
\end{equation}
where $\tau^{-\alpha}=\module{\tau}^{-\alpha}e^{-2i\pi\alpha\theta}$, $\log\tau=(\log\module{\tau})+\theta$, with $\theta\in{}]\epsilon-1,\epsilon[$ and $c_{\eta,\ell}^{(m)}\in\CC$. The coefficients $c_{\eta,\ell}^{(m)}$ only depend on the class $[\eta]$ of $\eta$ in $\gr_\alpha^VG$, so we will denote them by $c_{[\eta],\ell}^{(m)}$, and we have
\begin{equation}\label{eq:cetam}
c_{[\eta],\ell}^{(m)}=c_{(2i\pi N)^m[\eta],\ell}^{(0)}.
\end{equation}

We will now characterize the Jordan blocks $B_{k_0}$ in $H_\alpha$. Such a Jordan block is characterized by the constant value $s$ of $s_w(\cbbullet)$, so that it will be convenient to denote such a block by $B_{\alpha,s}$.

\begin{lemme}\label{lem:charB}
For $[\eta]\in H_\alpha$, we have $[\eta]\in B_{\alpha,s}$ if and only if, for any $\ell=0,\dots,\mu-1$ and any $j\geq0$, we have
\[
c_{(2i\pi N)^j[\eta],\ell}^{(0)}=\zeta^{\ell s}\sum_m c_{(2i\pi N)^{j+m}[\eta],0}^{(0)}\frac{(-\ell/\mu)^m}{m!}.
\]
\end{lemme}

\begin{proof}
For $\eta=\omega_k$, denote $c_{[\omega_k],\ell}^{(m)}=c_{[\omega'_k],\ell}^{(m)}$ by $c_{k,\ell}^{(m)}$. Then we have
\begin{equation}\label{eq:ckl}
\int_{\Delta_\ell(\tau)}e^{-\tau f}\omega_k=\tau^{-\sigma_w(k)}\Big(\Big[c_{k,\ell}^{(m_k)}\frac{(-\log\tau)^{m_k}}{m_k!}+\cdots+c_{k,\ell}^{(0)}\Big]+o(1)\Big),
\end{equation}
where $m_k+1$ denotes the order of $[\omega_k]$ in $\gr_{\sigma_w(k)}^VG$ with respect to $2i\pi N$. Remark now that, as $\omega_k=u^{a(k)}\omega_0$ and $\module{a(k)}=k$, we have
\[
\int_{\Delta_\ell(\zeta^{-\ell}\module{\tau})}e^{-\zeta^{-\ell}\module{\tau} f}\omega_k=\int_{\zeta^\ell\Delta}e^{-\zeta^{-\ell}\module{\tau} f}\omega_k= \zeta^{k\ell}\int_{\Delta}e^{-\module{\tau} f}\omega_k.
\]
Hence, we get
\begin{multline*}
c_{k,\ell}^{(m_k)}\frac{(-\log\module{\tau}+\ell/\mu)^{m_k}}{m_k!}+\cdots+c_{k,\ell}^{(0)}\\
=\zeta^{\ell s_w(k)} \Big[c_{k,0}^{(m_k)}\frac{(-\log\module{\tau})^{m_k}}{m_k!}+\cdots+c_{k,0}^{(0)}\Big],
\end{multline*}
and in particular
\begin{equation}\label{eq:ck}
c_{k,\ell}^{(0)}=\zeta^{\ell s_w(k)} \cdot\sum_{m=0}^{m_k}\frac{(-\ell/\mu)^m}{m!}\, c_{k,0}^{(m)}.
\end{equation}
Therefore, any element $[\eta]$ in $B_{\alpha,s}$ satisfies the equality of Lemma \ref{lem:charB} for $j=0$, hence for any $j$.

Conversely, remark first that, if $[\eta]$ is fixed, then the equality of Lemma \ref{lem:charB} for any $j\geq0$ is equivalent to
\[
\psi_{[\eta],\ell}(\zeta^{-\ell}\module{\tau})=\zeta^{\ell s}\psi_{[\eta],0}(\module{\tau}),
\]
where $\psi_{[\eta],\ell}$ is defined by \eqref{eq:psieta} (two polynomials are equal iff all the corresponding derivatives at $0$ are equal).

Write $[\eta]=\sum\lambda_k[\omega'_k]$ in $H_\alpha$, denote $m'_{[\eta]}=\max_{k\mid\lambda_k\neq0}m_k$ and put \hbox{$K_{[\eta]}=\{k\mid m_k=m'_{[\eta]}\}$}. Notice that, for $k,k'\in K_{[\eta]}$, we have $s_w(k)\neq s_w(k')$. If $m'_{[\eta]}>m_{[\eta]}$, we have $\sum_{k\in K_{[\eta]}}\zeta^{\ell s_w(k)}\lambda_k c_{k,0}^{(m_k)}=0$ for any $\ell=0,\dots,\mu-1$. It follows that $\lambda_k c_{k,0}^{(m_k)}=0$, hence $\lambda_k=0$, for any $k\in K_{[\eta]}$, a contradiction. Therefore $m'_{[\eta]}=m_{[\eta]}$. Argue similarly to show that $K_{[\eta]}$ is reduced to one element, denoted by $k_{[\eta]}$, and that $s_w(k_{[\eta]})=s$. Apply the lemma by induction on $m_{[\eta']}$ to $[\eta']=[\eta]-\lambda_{k_{[\eta]}}[\omega'_{k_{[\eta]}}]$.
\end{proof}

\subsection{Isomorphism between nearby cycles}
The multivalued cycles $\Delta_\ell(\tau)$ form a basis of the space of multivalued global sections of $\cH$, that we denote by $\psi_\tau\cH$. This basis defines the integral (hence the real) structure on $\psi_\tau\cH$.

Denote by $\cN_{\alpha,p}$ the space of linear combinations with meromorphic coefficients of germs at $\tau=0$ of the multivalued functions $e_{\alpha,q}=\tau^{\alpha}(-\log\tau)^q/q!$ ($q\leq p$). For $p$ large enough (here $p\geq n+1$ is enough), the map
\begin{align*}
\varphi:V_\alpha G&\to V_0(G\otimes\cN_{\alpha,p})\\
\eta&\mto \sum_{j=0}^p[2i\pi(\tau\partial_\tau+\alpha)]^j\eta\otimes e_{\alpha,j}
\end{align*}
induces an isomorphism
\[
\gr_\alpha^VG\isom \ker\big[\tau\partial_\tau:\gr_0^V(G\otimes\cN_{\alpha,p})\to\gr_0^V(G\otimes\cN_{\alpha,p})\big].
\]
As $G$ is regular at $\tau=0$, there exists a perturbation $\eta\mto\psi(\eta)\in V_{<0}(G\otimes\cN_{\alpha,p})$ such that $\varphi(\eta)+\psi(\eta)\in\ker\big[\tau\partial_\tau:G\otimes\cN_{\alpha,p}\to G\otimes\cN_{\alpha,p}\big]$.

Recall (see, \eg \cite{M-S86}) that $\ker\big[\tau\partial_\tau:G\otimes\cN_{\alpha,p}\to G\otimes\cN_{\alpha,p}\big]$ is identified with $H_\alpha$. Set $\cN=\oplus_{\alpha\in[0,1[}\cN_{\alpha,n+1}$. Given a section $\lambda$ of $H=\ker\big[\tau\partial_\tau:G\otimes\cN\to G\otimes\cN\big]$ and a section $\delta$ of $\psi_\tau\cH$, choose a representative $\wt \lambda$ of $\lambda$ in $\Omega^n(U)\otimes_\CC\cN$. Then $\int_\delta e^{-\tau f}\wt \lambda\in\CC$. Then (see Appendix) $\lambda$ belongs to $H_\QQ$ if and only if, for any $\ell=0,\dots,\mu$ and some nonzero $\tau$, we have
\begin{equation}\label{eq:comparaison}
\int_{\Delta_\ell(\tau)} e^{-\tau f}\wt\lambda\in\QQ.
\end{equation}
For $\eta\in V_\alpha G$ and $\wt\lambda=\varphi(\eta)+\psi(\eta)$, and using \eqref{eq:psieta}, one finds
\[
\int_{\Delta_\ell(\tau)} e^{-\tau f}\varphi(\eta)=c_{\eta,\ell}^{(0)}+ o(1).
\]
As a consequence, the conjugate $\ov{[\eta]}$ of $[\eta]$ satisfies
\[
c_{\ov{[\eta]},\ell}^{(0)}=\ov{c_{[\eta],\ell}^{(0)}}.
\]
It follows now from Lemma \ref{lem:charB} that
\[
\ov{B_{\alpha,s}}=
\begin{cases}
B_{1-\alpha,\mu-s}&\text{if }\alpha\in{}]0,1[,\\
B_{0,\mu-s}&\text{if }\alpha=0.
\end{cases}
\]
As $s_w(\ov k_0)=\mu-s_w(k_0)$, this ends the proof of Lemma \ref{lem:conjug}.\hfill\qed

\section*{Appendix}
\setcounter{section}{0}
\refstepcounter{section}
\renewcommand{\thesection}{\Alph{section}}

In this appendix, we explain with some details why the real structure on $H$ as defined by \eqref{eq:comparaison} is indeed the real structure used in \cite{Sabbah96a} to define the Hodge structure on $H$. We will need to recall some notation and results of \cite{Sabbah96a}.

We will denote by $U$ a smooth quasi-projective variety and by \hbox{$f:U\to\Afu$} a regular function on $U$. We denote by $t$ the coordinate on the affine line $\Afu$. We also fix an embedding $\kappa:U\hto \cX$ into a smooth projective variety such that there exists an algebraic map $F:\cX\to\PP^1$ extending~$f$. We have a commutative diagram, where the right part is Cartesian, thus defining $X$ as a fibred product,
\[
\xymatrix@C=1.5cm{
U\ar@{^{ (}->}[r]^-j\ar@/^2pc/[rr]^-{\kappa}\ar[dr]_-{f}&X\ar@{^{ (}->}[r]\ar[d]^-{F_{|X}}\ar@{}[dr]|{\scriptstyle\square}&\cX\ar[d]^-F\\
&\Afu\ar@{^{ (}->}[r]&\PP^1
}
\]

Denote by $\epsilon:\wt\PP^1\to\PP^1$ the real blow-up of $\PP^1$ centered at $\infty$ ($\wt\PP^1$ is diffeomorphic to a closed disc) and by $\wt F:\wt\cX\to\wt\PP^1$ the fibre-product of $F$ with the blowing-up $\epsilon$. Denote by $\wt\kappa$ the inclusion $U\hto\wt\cX$.

Denote by $S^1$ the inverse image of $\infty$ by the blowing-up $\epsilon$. Let $\Afuh$ be an affine line with coordinate $\tau$. Denote by $L^{\prime+}$ the closed set of $S^1\times\Afuh\subset\wt\PP^1\times\Afuh$ defined by $\reel(e^{i\theta}\tau)\geq0$, with $\theta=\arg t$ and where $t$ is the coordinate on $\Afu=\PP^1\moins\{\infty\}$, and set $L^-=\wt\PP^1\times\Afuh\moins L^{\prime+}$. For $\tau\neq0$, denote by $L^{\prime+}_\tau,L^-_\tau\subset\wt\PP^1$ the fibre of $L^{\prime+},L^-$ over $\tau$.

\begin{figure}[htb]
\centerline{\includegraphics[scale=.8]{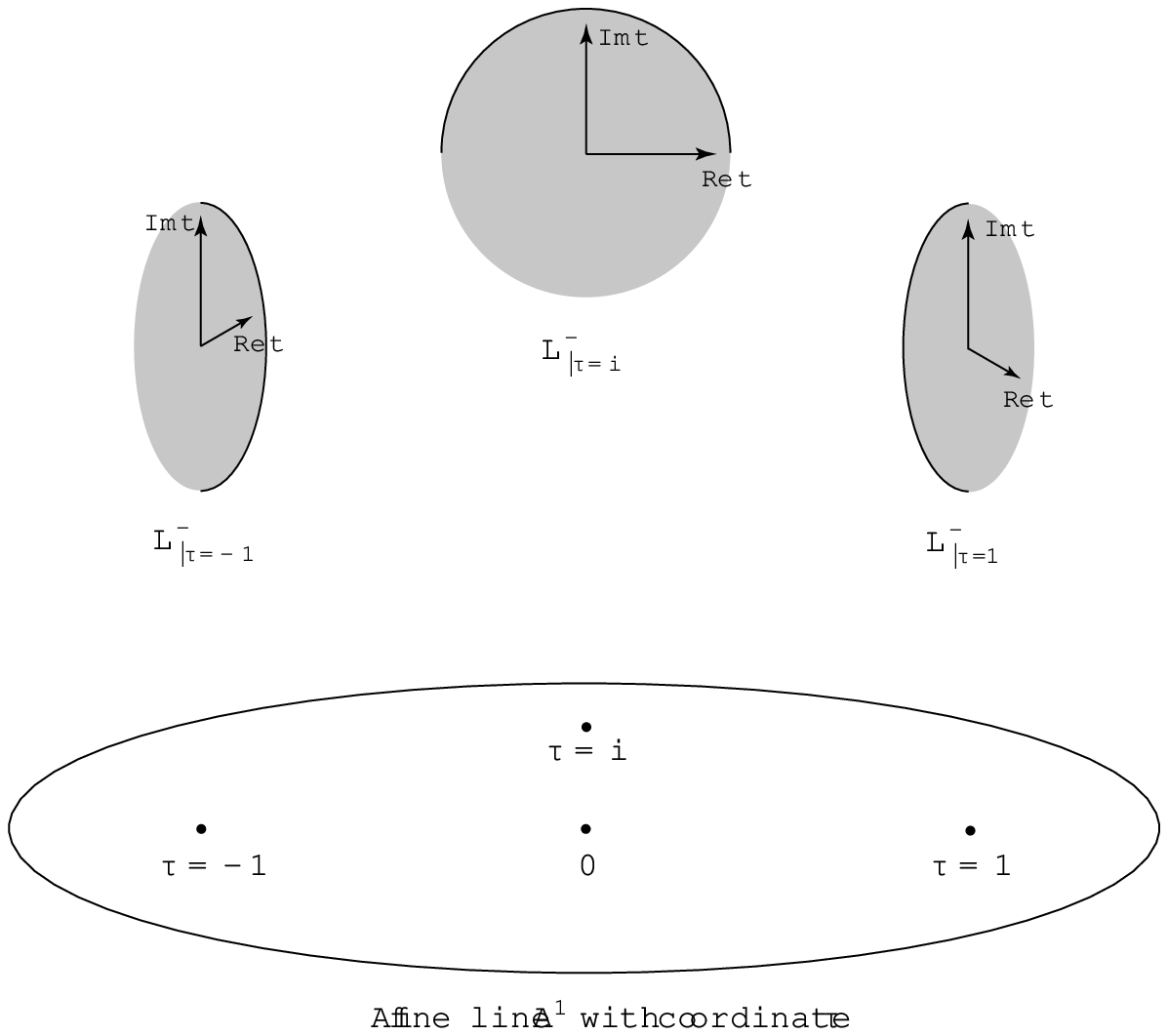}}
\end{figure}

Denote similarly by $L^{\prime+}_{\wt\cX},L^-_{\wt\cX}\subset\wt\cX\times \Afuh$ (\resp $L^{\prime+}_{\wt\cX,\tau},L^-_{\wt\cX,\tau}\subset\wt\cX$) the inverse image of the corresponding sets by $\wt F\times\id_{\Afuh}$ (\resp by $\wt F$).

We denote by $\alpha:\Afu\times\Afuh\hto L^-$ and $\beta:L^-\hto\wt\PP^1\times\Afuh$ (\resp $\alpha_\tau:\Afu\hto L^-_\tau$ and $\beta_\tau:L^-_\tau\hto\wt\PP^1$) the inclusions, and by the same letters the corresponding inclusions
\[
\alpha:X\times\Afuh\hto L^-_{\wt\cX}\quad\text{and}\quad \beta:L^-_{\wt\cX}\hto\wt\cX\times\Afuh,
\]
\resp
\[
\alpha_\tau:X\hto L^-_{\wt\cX,\tau}\quad\text{and}\quad \beta_\tau:L^-_{\wt\cX,\tau}\hto\wt\cX.
\]
Therefore we have $\beta_\tau\circ\alpha_\tau\circ j=\wt\kappa$.

In \cite[(1.8)]{Sabbah96a}, we have defined the Fourier transform $\gF_F(\bR j_*\CC_U)$ as the following complex on $\cX\times\Afuh$ (there is a shift by $1$ in \loccit., that we do not introduce here):
\[
\gF_F(\bR j_*\CC_U)\defin \bR\epsilon_*\,\beta_!\,\bR\alpha_*\,\bR j_*\CC_{U\times\Afuh},
\]
where we still denote by $j$ (\resp $\kappa$) the inclusion $U\times\Afuh\hto X\times\Afuh$ (\resp $U\times\Afuh\hto \cX\times\Afuh$). This complex has a natural $\QQ$-structure (replace $\CC_U$ with $\QQ_U$). This induces a $\QQ$-structure on the nearby cycle complex $\psi_\tau\gF_F(\bR j_*\CC_U)=\psi_\tau\gF_F(\bR j_*\QQ_U)\otimes_\QQ\CC$.

Denote by $\cE^{-\tau f}$ the algebraic $\cD_{U\times\Afuh}$-module $\cO_{U\times\Afuh}e^{-\tau f}$ (\ie the $\cO_{U\times\Afuh}$-module $\cO_{U\times\Afuh}$ with connection $e^{\tau f}\circ d\circ e^{-\tau f}$). The quasi-isomorphism
\begin{equation}\label{eq:Sabbah2.2}
\DR^\an_{\wt\cX\times\Afuh}(\kappa_+\cE^{-\tau f})\isom\gF_F(\bR j_*\CC_U)
\end{equation}
constructed in \cite[Th\ptbl2.2]{Sabbah96a} is then used to define the $\QQ$-structure on the complex (of sheaves on $\cX$) $\psi_\tau\DR^\an_{\wt\cX\times\Afuh}(\kappa_+\cE^{-\tau f})$ [on the other hand one uses the $V$-filtration relative to $\tau=0$ on $\kappa_+\cE^{-\tau f}$ to construct the Hodge filtration on this complex]. By $\DR$ we mean the usual de~Rham complex, starting in degree $0$.

Denote by $\wt\cO_{\Afuh}^\an$ the sheaf of multivalued holomorphic functions on $\Afuh\moins\{0\}$. Then
\[
\psi_\tau\DR^\an_{\wt\cX\times\Afuh}(\kappa_+\cE^{-\tau f})=i_{\tau=0}^{-1}\DR^\an_{\wt\cX\times\Afuh}(\kappa_+\cE^{-\tau f}\otimes_{\wh p^{-1}{\cO_{\Afuh}^\an}}\wt\cO_{\Afuh}^\an),
\]
where $\wh p:\cX\times\Afuh\to\Afuh$ denotes the projection and $i_{\tau=0}:\cX\times\{0\}\hto\cX\times\Afuh$ denotes the inclusion (see, \eg \cite[(4.9.4)]{M-S86}). We are interested in analyzing the $\QQ$-structure on the cohomology of $\bR\Gamma\big(\cX,\psi_\tau\DR^\an_{\wt\cX\times\Afuh}(\kappa_+\cE^{-\tau f})\big)$. Use $C^\infty$ forms on $\cX$ to identify it with
\[
\Gamma\big(\cX,i_{\tau=0}^{-1}\cE_{\cX\times\Afuh}^{\cbbullet}(\kappa_+\cE^{-\tau f}\otimes_{\wh p^{-1}{\cO_{\Afuh}^\an}}\wt\cO_{\Afuh}^\an)\big),
\]
with $n=\dim \cX$. Similarly, denote by $\wt\CC_{\Afuh}$ the sheaf of multivalued local sections of $\CC_{\Afuh}$ (\ie local sections on the universal covering of $\Afuh\moins\{0\}$). Then $\psi_\tau\gF_F(\bR j_*\CC_U)$ is equal to $i_{\tau=0}^{-1}\big(\gF_F(\bR j_*\CC_U)\otimes\wh p^{-1}\wt\CC_{\Afuh}\big)$.

In order to know that the cohomology class of a closed multivalued section of $\wh p_*\cE_{\cX\times\Afuh}^{\cbbullet}(\kappa_+\cE^{-\tau f}\otimes_{\wh p^{-1}{\cO_{\Afuh}^\an}}\wt\cO_{\Afuh}^\an)$ is rational, one has to compute its image in $\bR\wh p_*\gF_F(\bR j_*\CC_U)\otimes\wt\CC_{\Afuh}$ and decide whether its class is rational or not. As the section is closed, it is enough to verify this after restricting to some (or any) $\tau\neq0$. Therefore, we need to compute the map \eqref{eq:Sabbah2.2} after restricting to some fixed nonzero $\tau$. In \eqref{eq:comparaison}, we apply this computation to the multivalued form~$e^{-\tau f}\wt\lambda$.

Denote by $\cE_{\wt\cX}$ the sheaf of $C^\infty$ functions (in the sense of Whitney) on $\wt\cX$, by $\cE^{\mod}_{\wt\cX}$ the sheaf on $\wt\cX$ of $C^\infty$ functions on $U$ which have moderate growth along $\wt\cX\moins U$, and by $\cE^{\mod,-_\tau}_{\wt\cX}$ the subsheaf of functions which moreover are infinitely flat along $L^{\prime+}_{\wt\cX,\tau}$.

On the other hand, denote by $\cC^\cbbullet_{U\cup L^{\prime+}_{\wt\cX,\tau},L^{\prime+}_{\wt\cX,\tau}}$ the complex of sheaves on $\wt\cX$, consisting of germs on $\wt\cX$ of relative singular cochains (\ie germs of singular cochains in $U\cup L^{\prime+}_{\wt\cX,\tau}$ with boundary in $L^{\prime+}_{\wt\cX,\tau}$).

By the de~Rham theorem, the integration of forms induces a quasi-isomorphism of complexes $\int:\cE^\cbbullet_U\to\cC_U^\cbbullet\otimes_\ZZ\CC$; moreover, the natural morphism $\cE^{\mod,\cbbullet}_{\wt\cX}\to(\alpha_\tau\circ j)_*\cE^\cbbullet_U$ is a quasi-isomorphism, so the integration morphism $\int:\cE^{\mod,\cbbullet}_{\wt\cX}\to(\alpha_\tau\circ j)_*\cC_U^\cbbullet\otimes_\ZZ\CC$, which is obtained by composing both morphisms, is a quasi-isomorphism.

Similarly, we have a commutative diagram
\[
\xymatrix{
\beta_{\tau,!}\cE^{\mod,\cbbullet}_{\wt\cX}\ar[r]^-\sim\ar[d]_-\wr& \cE^{\mod,-_\tau,\cbbullet}_{\wt\cX}\ar[dr]^-\int\\
\beta_{\tau,!}(\alpha_\tau\circ j)_*\cE^\cbbullet_U\ar[r]_-\int^-\sim& \beta_{\tau,!}(\alpha_\tau\circ j)_*\cC_U^\cbbullet\otimes_\ZZ\CC\ar[r]^-\sim &\cC^\cbbullet_{U\cup L^{\prime+}_{\wt\cX,\tau},L^{\prime+}_{\wt\cX,\tau}} \otimes_\ZZ\CC
}
\]
Hence we get:

\begin{proposition}[A variant of the de~Rham theorem]
Both complexes $\cC^\cbbullet_{U\cup L^{\prime+}_{\wt\cX,\tau},L^{\prime+}_{\wt\cX,\tau}}\otimes_\ZZ\CC$ and $\cE^{\mod,-_\tau,\cbbullet}_{\wt\cX}$ are quasi-isomorphic to $\beta_{\tau,!}\bR\alpha_{\tau,*}\bR j_*\CC_U$. Moreover, the integration of forms induces a natural quasi-isomorphism of complexes
\[
\int:\cE^{\mod,-_\tau,\cbbullet}_{\wt\cX}\isom\cC^\cbbullet_{U\cup L^{\prime+}_{\wt\cX,\tau},L^{\prime+}_{\wt\cX,\tau}}\otimes_\ZZ\CC.\qedhere
\]
\end{proposition}

Now, given a section of $\cE^{\mod,\cbbullet}_{\wt\cX}\otimes j_+\cE^{-\tau f}$, \ie a section of $\cE^{\mod,\cbbullet}_{\wt\cX}$ multiplied by $e^{-\tau f}$, it is also a section of $\cE^{\mod,-_\tau,\cbbullet}_{\wt\cX}$, and its image by \eqref{eq:Sabbah2.2} is nothing but its integral, according to the previous commutative diagram.

\backmatter
\providecommand{\bysame}{\leavevmode ---\ }
\providecommand{\og}{``}
\providecommand{\fg}{''}
\providecommand{\smfandname}{\&}
\providecommand{\smfedsname}{\'eds.}
\providecommand{\smfedname}{\'ed.}
\providecommand{\smfmastersthesisname}{M\'emoire}
\providecommand{\smfphdthesisname}{Th\`ese}

\end{document}